\magnification=\magstep 1
\hbadness=10000
\hsize=30truecc
\parskip=0pt
\baselineskip=18 pt
\hsize=6.5truein
\vsize=8.5truein
\hfuzz=2pt
\overfullrule=0pt   
\output={\plainoutput}
\def\abstract{\vskip 0.5truecm
\centerline{\bf Abstract}\vskip 5pt}

\def\section#1.#2.{\baselineskip=12pt\goodbreak\vskip1truecm\hangindent
            \hangafter=1\noindent{\bf #1.}\ \ {\bf #2.}\nobreak\baselineskip=20pt\ignorespaces}
\headline={\hss{}\hss}
\footline={\ifnum\pageno<2 \hss{}\hss\else\tenrm\hss\folio\hss\fi}

\newcount\refno
\refno=0
\def\ref{\advance\refno by 1\medbreak\item{\the\refno.}}

\def\bibs{\vskip 1truecm\baselineskip12pt\parskip 3 pt
                         \centerline{\bf References}\par\nobreak}

\def\noi{\noindent}
\def\squar{\vbox{\hrule\hbox{\vrule height 6pt \hskip
6pt\vrule}\hrule}}
\def\sqr{{\unskip\nobreak\hfil\penalty50\hskip2em\hbox{}\nobreak\hfil
{\squar}\parfillskip=0pt\finalhyphendemerits=0\par}}

\def\Thm #1 {\medbreak\noi{\bf Theorem #1}\ \ }
\def\Lemma#1 {\medbreak\noi{\bf Lemma #1}\ \ }
\def\Cor#1 {\medbreak\noi{\bf Corollary #1}\ \ }
\def\Prop#1 {\medbreak\noi{\bf Proposition #1}\ \ }

\def\Remarks#1 {\medbreak\noi{\bf Remarks #1}\ \ }

\def\Def #1 {\medbreak\noi{\bf Definition #1}\ \ }
\def\proof{\noi{\bf Proof}:\ \ }
\def\endproof{\sqr\par}
\def\der#1.{^{(#1)}}

\def\sig{\sigma}
\def\eps{\epsilon}
\def\ga{\gamma}

\def\a{\alpha}
\def\b{\beta}
\def\th{\theta}

\def\B{{\cal{B}}}
\def\F{{\cal{F}}}
\def\R{{\cal R}}
\def\N{{{\cal N}}}

\def\dover{\buildrel d \over =}

\def\berg{BerG$\mr$\ }
\def\bergp{BerG$\mrp$\ }

\def\Xtone{X_{t-1}}
\def\Xtk{X_{t-k}}
\def\Xt{\{X_t\}}
\def\epst{\{\eps_t\}}

\def\bzp{$\bf Z_+$}
\def\mr{(m,r)}
\def\mrp{(m',r')}
\def\mrone{(m_1,r_1)}
\def\mrtwo{(m_2,r_2)}
\def\mrastk{(m,r)^{*k}}
\def\mrastn{(m,r)^{*n}}
\def\MR{(M,R)}
\def\MRSq{[M,R]}

\def\psimr{\psi_{m,r}}
\def\psimrp{\psi_{m',r'}}
\def\psimrone{\psi_{m_1,r_1}}
\def\psimrtwo{\psi_{m_2,r_2}}
\def\psione{\psi_{1,r}}
\def\psimrk{\psi_{m,r}^{(k)}}

\def\pkn{p_k^{(n)}}

\def\finar1{$F$-INAR$\,{(1)}$\ }
\def\inar1{INAR$\,(1)$}
\def\nginar1{NGINAR$\,(1)$}

{\bf \centerline{ Expectation thinning operators based on linear fractional}
\centerline{ probability generating functions}}

\vskip 1 cm

{\baselineskip =12pt\halign to \hsize{\qquad#\hfil &\qquad#\hfil\cr
Emad-Eldin A. A. Aly & Nadjib Bouzar\cr
Department of Statistics and O.R. & Department of Mathematical Sciences\cr
Kuwait University & University of Indianapolis\cr
P.O.B. 5969, Safat 13060 & Indianapolis, IN 46227\cr
Kuwait & U.S.A.\cr
\underbar{Email}: emad@kuc01.kuniv.edu.kw & \underbar{Email}: nbouzar@uindy.edu\cr}}

\vskip1.5cm \abstract We introduce a two-parameter expectation thinning operator based on a linear fractional probability generating function. The operator is then used to define a first-order integer-valued autoregressive \inar1 process. Distributional properties of the \inar1 process are described. We revisit the Bernoulli-geometric \inar1 process of Bourguignon and Wei{\ss} (2017) and we introduce a new stationary \inar1 process with a compound negative binomial distribution. Lastly, we show how a proper randomization of our operator leads to a generalized notion of monotonicity for distributions on \bzp.  

{\baselineskip 12 pt \footnote{}{\hskip -.3 in {\it Key words and phrases\/}: semigroup, \inar1 process, re-parameterization, stationarity, monotonicity.}}

\vfill\eject

\noi\nobreak {\bf 1. Introduction}
\smallskip
Thinning operators have been successfully used in the last thirty years to model time series for count data. These operators preserve the discrete nature of the variates and play the role of a generalized multiplication in the equations that govern integer-valued autoregressive moving average (INARMA) models. 

Historically, the binomial thinning operator $\otimes_\B$ of Steutel and van Harn (1979) was the first operator used to construct thinning-based INARMA models. It is defined as follows.

\Def 1.1. Let $\a\in (0,1)$ and $X$ a \bzp-valued random variable. Then
$$
\a\otimes_\B X=\sum_{i=1}^X B_i,\leqno(1.1)
$$
where $(B_i, i\ge 1)$ is a sequence of iid Bernoulli($\a$) random variables, independent of $X$.

As noted in Wei{\ss} (2008), binomial thinning-based INARMA models perform well with Poissonian count data, but not as well with variates that exhibit overdispersion or underdispersion. We refer the reader to the excellent survey articles by McKenzie (2003), Wei{\ss} (2008), and Scotto et al. (2015) for a deeper discussion of these issues.  

Alternatives to the binomial thinning operators were proposed by several authors. These generalized thinning operators have been  designed to deal with count data that show overdispersion or underdispersion due in particular to a deflation or an inflation of zeros. We will follow Zhu and Joe (2003) and refer to these operators as expectation thinning operators in the sense that at any given time, the action of the operator on a variate yields a smaller expected count than the value of the variate at that time.

The focus of this article will be on the expectation thinning operators based on linear fractional probability generating functions (pgf's). These operators have been particularly useful in modeling stationary first order integer-valued autoregressive (\inar1) processes with geometric, negative binomial, and Poisson-geometric marginal distributions.

In Section 2, we establish that any nondegenerate linear fractional pgf $f(s)$ gives rise, via a suitable re-parameterization, to a two-parameter operator that enjoys a useful semigroup property as well as the standard linearity properties for the conditional means and variances of variates. Moreover, the operator will be of the expectation thinning type if $0<f'(1)<1$. We show that several expectation thinning operators based on specific linear fractional pgf's arise as special cases of our operator (via re-parameterization). These operators are individually referenced at the end of the section. 
   
In Section 3, we use the thinning version of our operator to define a first-order integer-valued autoregressive (\inar1) process. We state the main distributional properties of the process. We revisit the Bernoulli-geometric \inar1 process of Bourguignon and Wei{\ss} (2017) and show that the range of admissible values of its parameters extends to a larger set. We also propose a stationary \inar1 model with the zero-modified marginal distribution of Barreto-Souza (2015). Lastly, we introduce a new stationary \inar1 process with a compound negative binomial distribution and derive the distribution of its innovation sequence. 

In Section 4, we show how a proper randomization of our operator leads to a generalized notion of monotonicity for distributions on \bzp. Our results are to be seen as generalizations of $\a$-monotonicity introduced by Steutel (1988) (based on binomial thinning) and of $(\rho,\a)$-generalized mononoticity of Jazi and Alamatsaz (2012) (based on an expectation thinning operator driven by a linear fractional pgf). 

\vskip .5 cm

\noi\nobreak {\bf 2. A two-parameter expectation thinning operator}
\smallskip

Let $\displaystyle f(s)={a+bs\over c+ds}$, $s\in [0,1]$, be a linear fractional pgf, with $f(0)<1$. A straightforward power series argument shows that $f(s)$, relabeled henceforth as $\psimr(s)$, can be rewritten  in the form
$$
\psimr(s)=1-m{1-s\over 1+r(1-s)}\quad s\in [0,1], \leqno(2.1)
$$
where $m=f'(1)$, $r\ge 0$, and $0< m\le r+1$.

Let 
$$
\R=\{\mr\in {\bf R}^2: r\ge 0 \hbox{\ and \ } 0<m\le r+1\}.\leqno(2.2)
$$

We recall that a \bzp-valued random variable $X$ is said to have a T-geometric$(p)$ distribution, $p\in (0,1)$ (and T for truncated at zero), if its probability mass function (pmf) is $P(X=k)=p(1-p)^{k-1}$, $k\ge 1$.

We start out by listing several useful properties of the pgf $\psimr(s)$. 

\Prop 2.1. Let $\mr\in\R$ and $Z$ a \bzp-valued random variable with pgf $\psimr(s)$. 

\noi(i) The pmf of $Z$ is
$$
p_k=\cases{
1-{m\over 1+r}\quad & if $k=0$\cr
\strut
{mr^{k-1}\over (1+r)^{k+1}}\quad & if $k\ge 1$.\cr
}
\leqno(2.3)
$$
\noi (ii) $Z$ admits the representations
$$
Z\dover BW\dover \sum_{i=1}^{B'} W_i',\leqno(2.4)
$$
where $B$ and $W$ are independent, $\{W_i\}$ is a sequence of iid random variables independent of $B'$, $B$ and $B'$ are Bernoulli$\bigl({m\over r+1}\bigr)$, $W$ and the $W_i$'s are T-geometric$\bigl({1\over r+1}\bigr)$.

\noi (iii) The mean and variance of $Z$ are
$$
E(Z)=m\quad  \hbox{and} \quad Var(Z)=m(2r+1-m),
$$
and the dispersion index $I_Z={Var(z)\over E(Z)}=2r+1-m$ indicates equidispersion of $\{p_k\}$ if $m=2r$, underdispersion if $m>2r$ and overdispersion if $m<2r$.
  
\noi (iv) Let $n\ge 1$. The pmf of the $n$-fold convolution of $\{p_k\}$ of (2.3) is
$$
\pkn=\Bigl(1-{m\over r+1}\Bigr)^n\Bigl({r\over r+1}\Bigr)^k\sum_{i=0}^{\min(k,n)}{n\choose i}{k-1\choose i-1}\Bigl({m\over r(r+1-m)}\Bigr)^i.\leqno(2.5)
$$
\proof The proof of (i)-(iii) is a simple exercise. For (iv), we note that if $Y$ is a \bzp-valued random variable with pgf $\psimr^n(s)$, then it admits the representation $Y\dover\sum_{i=1}^N Y_i$, where $N\sim \hbox{Binomial}\bigl(n,m/(r+1)\bigr)$ and $\{Y_i\}$ is a sequence of iid random variables, independent of $N$, and such that $Y_i\sim \hbox{T-Geometric}\bigl(1/(r+1)\bigr)$. A standard conditioning argument leads to (2.5).\endproof

The pmf  (2.1) for $r<m\le r+1$ appears in Bourguignon and Wei{\ss} (2017) under a different parameterization (see additional details at the end of the section). The authors named it the BerG distribution as it results from the convolution of a Bernoulli $(m-r)$ distribution and a (non-truncated at zero) geometric $\bigl({1\over r+1}\bigr)$ distribution. We extend the label to any $\mr\in\R$  and will refer to a \berg distribution as the distribution with pmf (2.1) (or pgf $\psimr(s)$).

Next, we define a binary operation on $\R$ as follows:
$$
\mr \ast \mrp=(mm', r+r'm)\qquad \mr, \mrp\in\R.\leqno(2.6)
$$ 
$\R$ equipped with the operation $(\ast)$ is a semigroup. Indeed, $\R$ is closed under $(\ast)$ as $mm'\le mr'+m\le mr'+r+1$. It is easily seen that $(\ast)$ is associative and that it admits $(1,0)$ as its neutral element. In general, $(\ast)$ is not commutative, In fact, if $\mr$ and $\mrp$ are in $\R$, then $\mr \ast \mrp=\mrp \ast \mr$ if and only if $r(1-m')=r'(1-m)$. We note that $(\ast)$ is commutative when restricted to the following sub-semigroups of $\R$: $A=\{\mr\in \R: 0<m\le 1 \hbox{\ and\ }r=0\}$, $B=\{\mr\in \R: m=1\}$, $C=\{\mr\in \R: m=r+1\}$.

Let $\mr\in \R$. We define
$$
\mrastk=\underbrace{\mr\ast\mr\ast\cdots\ast\mr}_{k \ \hbox{times}} \qquad (k\ge 1),\leqno(2.7)
$$
with $\mr^{\ast 0}=(1,0)$.

Assume $\mr\in \R$. By (2.6) and a simple induction argument, we have
$$
\mrastk=\bigl(m^k, r s_k\bigr), \qquad s_k=\sum_{j=0}^{k-1} m^j \quad (k\ge 1).\leqno(2.8)
$$

The family of pgf's $\Psi=(\psimr(\cdot), \mr\in \R)$ enjoys the following semigroup property (proof is omitted). 

\Prop 2.2. For any $\mr$ and $\mrp$ in $\R$,
$$
\psimrp(\psimr(s))=\psi_{\mr\ast \mrp}(s).\leqno(2.9)
$$

We define the iterates of $\psimr(s)$, $\mr\in \R$, by
$$
\psimrk(x)\cases{\psimr(s) & if $k=1$\cr
\psimr((\psimr^{(k-1)}(s)) & if $k\ge 2$.\cr 
}\leqno(2.10)
$$ 
We deduce by (2.7)-(2.9) and an induction argument that 
$$
\psimrk(s)=\psi_{\mrastk}(s)\quad (k\ge 1),\leqno(2.11)
$$ 

We now introduce a two-parameter operator that acts on \bzp-valued random variables.

\Def 2.3. Let $X$ be a \bzp-valued random variable and $\mr\in \R$. Then
$$
\mr\odot X=\sum_{i=1}^X Z_i,\leqno(2.12)
$$
where $\{Z_i\}$ is a sequence of iid \bzp-valued random variables independent of $X$ and with marginal pgf $\psimr(s)$ of (2.1). If $0<m<1$, we will refer to $\odot$ as an {\it expectation thinning operator}.

If $Q(s)$ is the pgf of $X$, then the pgf $P(s)$ of $\mr\odot X$ satisfies 
$$
P(s)=Q(\psimr(s)).\leqno(2.13)
$$

The operator $\odot$ enjoys the following closure property.

\Prop 2.4. Let $X$ be a \bzp-valued random variable and let $\mr$ and $\mrp$ be in $\R$. Then
$$
\mr\odot\bigl(\mrp\odot X\bigr)\dover (\mr\ast\mrp)\odot X.\leqno(2.14)
$$
\proof By (2.13), the pgf $\phi(s)$ of the left-hand side of (2.14) satisfies $\phi(s)=Q(\psimrp(\psimr(s)))$, where $Q$ is the pgf of $X$. It follows by (2.9) that  is $\phi(s)=\phi(\psi_{m'',r''}(s))$, with $(m'',r'')=\mr\ast\mrp$.\endproof

Let $X$ be a \bzp-valued random variable and $\mr\in \R$. We define the $k$-fold action of $\mr\odot (\cdot)$ on $X$ by
$$
Y_k=\cases{
\mr\odot X \quad & if $k=1$\cr
\mr\odot Y_{k-1}\quad & if $k\ge 2$\cr
}
$$
We will use the notation below without further reference:
$$
Y_k=\underbrace{\mr\odot \mr\odot \cdots\odot \mr\odot}_{k\ times}X.
$$  
Proposition 2.4 and an induction argument lead to the following result.

\Cor 2.5. Let $X$ be a \bzp-valued random variable and $\mr\in\R$. Then
$$
\mr\odot \mr\odot \cdots\odot \mr\odot X\dover \mrastk\odot X,\leqno(2.15)
$$
where $\mrastk$ is as in (2.8).

We note that the expectation thinning operator $(m,0)\odot X$, $0<m <1$, becomes the binomial thinning operator $m\otimes_\B X$ of (1.1) as in this case the $Z_i$'s in (2.12) will have a common Bernoulli($m$) distribution.

For $m=1$ and $r>0$ the $\odot$ operator of (2.12) becomes a special case of the van Harn et al. (1982) $\odot_\F$ operator, where $\F=(F_r(\cdot), r\ge 0)$ is a continuous semigroup of pgf's. Indeed, we see by (2.9) that $\Psi_1=(\psione(\cdot), r\ge 0)$ forms a continuous semigroup of pgf's. In this case the van Harn et al. operator, which we denote by $\otimes_{\Psi_1}$, is defined by
$$
e^{-r}\otimes_{\Psi_1} X=\sum_{i=1}^X Z_i \qquad (r\ge 0),\leqno(2.16)
$$
where $X$ is a \bzp-valued random variable and $\{Z_i\}$ is a sequence of iid \bzp-valued random variables independent of $X$ and with marginal pgf $\psione(s)$. Since the pgf of $e^{-r}\otimes_\Psi X$ is $Q(\psione(s))$, where $Q$ is the pgf of $X$, we can conclude from (2.13) that 
$$
(1,r)\odot X \dover e^{-r}\otimes_{\Psi_1} X. \leqno(2.17)
$$

The operator $\odot$ becomes a single parameter operator if $m=r+1$ or $m=r$, with $Z_i\sim \hbox{T-Geometric}\bigl(1/(r+1)\bigr)$ when $m=r+1$ and $Z_i\sim \hbox{Geometric}\bigl(1/(r+1)\bigr)$ when $m=r$.
 
Noting that $(m,r)=(1,r)\ast (m,0)=(m,0)\ast(1,{r\over m})$, we obtain the following representations of the expectation thinning operator $\odot$ in terms of the operators $\otimes_{\Psi_1}$ and $\otimes_\B$.
 
\Prop 2.6. Let $X$ be a \bzp-valued random variable, $0<m\le 1$ and $r\ge 0$. Then
$$
\mr\odot X\dover e^{-r}\otimes_{\Psi_1}(m\otimes_\B X)\dover m\otimes_\B(e^{-{r\over m}}\otimes_{\Psi_1}X).\leqno(2.18) 
$$

We gather several properties of the operator $\odot$ in the following proposition. The proofs are omitted as they follow fairly straightforwardly from Proposition 2.1, equation (2.13), along with standard conditioning and pgf arguments for random summations.

\Prop 2.7. Let $\mr\in\R$ and $X$ a \bzp-valued random variable.

\noi (i) $E(\mr\odot X|X)=mX$.

\noi (ii) $E([\mr\odot X]^2|X)=(2r+1)mX+m^2X(X-1)$.

\noi (iii) $Var(\mr\odot X|X)=m(2r+1-m)X$.

\noi (iv) $E(\mr\odot X)=mE(X)$ and $Var(\mr\odot X)=m^2 Var(X)+m(2r+1-m)E(X)$.

\noi (v) For $k\ge 0$,
$$
P(\mr\odot X=k|X)=\cases
{\bigl(1-{m\over r+1}\bigr)^X \quad & if $k=0$\cr
\strut
\bigl(1-{m\over r+1}\bigr)^X \bigl({r\over r+1}\bigr)^k\sum_{i=0}^{\min(k,X)}{X\choose i}{k-1\choose i-1}\bigl({m\over r(r+1-m)}\bigr)^i\quad & if $k\ge 1$\cr
}.
$$

\noi (vi) If $Y$ is a \bzp-valued random variable independent of $X$, then
$$
\mr\odot(X+Y)\dover\mr\odot X + \mr\odot Y.
$$

We conclude the section by giving a fairly exhaustive list of expectation thinning operators based on a linear fractional pgf that appeared in the literature. We offer brief comments on how they relate to the $\odot$ operator.

(i) The binomial thinning operator $\otimes_\B$ of (1.1) is based on the pgf $\psi_{m,0}(s)=1-m+ms$, $0<m\le 1$, and as noted above, $m\otimes_\B(\cdot)\dover(m,0)\odot (\cdot)$.

(ii) The expectation thinning operator $A_{\a,\th}\circ (\cdot)$ of Aly and Bouzar (1994a, 1994b), $\a\in (0,1)$ and $\th\in [0.1]$, based on the pgf $\displaystyle \varphi_{\a,\th}(s) =1-\a{1-s\over 1-\th(1-\a)s}$: we have via re-parameterization $\varphi_{\a,\th}(s)=\psimr(s)$, with
$$
A_{\a,\th}\circ (\cdot)\dover \mr\odot (\cdot) \qquad m={\a\over 1-\th(1-\a)} \hbox{\ and\ } r={\th(1-\a)\over 1-\th (1-\a)}.\leqno(2.19)
$$

(iii) The expectation thinning operator $K(\a)\circ(\cdot)$ of Zhu and Joe (2003) based on the pgf $\varphi_{\a,\ga}(s)={(1-\a)+(\a-\ga)s\over (1-\a\ga)-(1-\a)\ga s}$, $\a\in [0,1]$ and $\ga\in [0,1)$: we have via re-parameterization $\varphi_{\a,\ga}(s)=\psimr(s)$, with
$$
K(\a)\circ(\cdot)\dover \mr\odot(\cdot) \qquad m=\a \hbox{\ and\ } r={\ga(1-\a)\over 1-\ga}.\leqno(2.20)
$$

(iv) The iterated thinning operator $\rho\star_\a(\cdot)$ of Wei{\ss} (2008) and Al-Osh and Aly (1992) based on the pgf $\varphi_{\a,\rho}(s)=1-{\a\rho(1-s)\over 1+\a-s}$, $0<\a,\ \rho<1$: we have via re-parameterization $\varphi_{\a,\rho}(s)=\psimr(s)$, with
$$
\rho\star_\a(\cdot)\dover\mr\odot(\cdot) \qquad m=\rho \hbox{\ and\ } r={1\over\a}.\leqno(2.21)
$$

(v) The negative binomial thinning operator $\a\odot_{\N\B}(\cdot)$ of Ristic et al. (2009) based on the pgf $\varphi_\a(s)={1\over 1+\a(1-s)}$, $\a\in [0,1)$: we have via re-parameterization $\varphi_\a(s)=\psimr(s)$, with
$$
\a\odot_{\N\B}(\cdot)\dover \mr\odot(\cdot) \qquad m=r=\a .\leqno(2.22)
$$

(vi) The operator $\pi \otimes_\rho(\cdot)$ of Jazi and Alamatsaz (2012) based on the pgf $\varphi_{\pi,\rho}(s)=1-{\pi(1-s)\over 1-\rho s}$, $\pi,\, \rho\in [0,1]$: we have via re-parameterization $\varphi_{\pi,\rho}(s)=\psimr(s)$, with
$$
\pi \otimes_\rho(\cdot)\dover \mr\odot(\cdot) \qquad m={\pi\over 1-\rho} \hbox{\ and\ } r={\rho\over1-\rho}.\leqno(2.23)
$$
The additional assumption $\pi+\rho<1$ makes $\otimes_\rho$ an expectation thinning operator.

(vii) The $\rho$-binomial thinning operator $\a \odot_{\rho\,\B}(\cdot)$ of Borges et al. (2016) based on the pgf $\varphi_{\a,\rho}(s)={1-(1-s)[\a(1+\rho)-\rho]\over 1+\rho(1-s)}$, $\rho\in [0,1)$ and $0\le \a<{1\over 1+\rho}$: we have via re-parameterization $\varphi_{\a,\rho}(s)=\psimr(s)$, with
$$
 \a \odot_{\rho\,\B}(\cdot)\dover \mr\odot(\cdot) \qquad m=\a(1+\rho) \hbox{\ and\ } r=\rho.\leqno(2.24)
$$

(viii) The $\rho$-negative binomial operator $\a \odot_{\rho\,\N\B}(\cdot)$ of Borges et al. (2017) based on the pgf $\varphi_{\a,\rho}(s)={1-\rho s\over 1-\rho s+\a(1-s)}$, $\rho\in [0,1)$, and $0< \a<1-\rho$: we have via re-parameterization $\varphi_{\a,\rho}(s)=\psimr(s)$, with
$$
\a \odot_{\rho\,\N\B}(\cdot)\dover \mr\odot(\cdot) \qquad m={\a \over 1-\rho} \hbox{\ and\ } r={\a+\rho\over 1-\rho}.\leqno(2.25)
$$

(ix) The Bourguignon and Wei{\ss} (2017) two-parameter operator $(\a,\b)\otimes(\cdot)$ based on the pgf $\varphi_{\a,\b}(s)={1-\a(1-s)\over 1+\b(1-s)}$, $\a\in (0,1)$, and $\b >0$: we have via re-parameterization $\varphi_{\a,\b}(s)=\psimr(s)$,
$$
(\a,\b)\otimes(\cdot)\dover \mr\odot(\cdot) \qquad m=\a+\b \hbox{\ and\ } r=\b.\leqno(2.26)
$$
Their operator is of the thinning type under the additional assumption $\a+\b<1$.

The thinning versions of the operators in (i)-(v) and (vii)-(ix) were primarily used to construct \inar1 processes with geometric, negative binomial, and Poisson-geometric marginal distributions. We refer to the original articles for more details. Jazi and Alamatsaz (2012) used their operator ((vi) above) to introduce a generalized notion of monotonicity for distributions on \bzp (more on this in Section 4).
   
Finally, we note that the linear fractional pgf's of the expectation thinning versions of the operators (i) and (iii)-(ix) can be  written (via suitable re-parameterizations) as $\varphi_{\a,\th}(s)$, the pgf of the operator $A_{\a,\th}$ of Aly and Bouzar (1994a, 1994b), for some $\a\in (0,1)$ and $\th\in [0.1]$, with the converse holding true only for operators (iii) and (vi).

\vskip .5 cm

\goodbreak\noi {\bf 3. An \inar1 process}
 
Let $\R_1=\{\mr\in \R: 0<m< 1 \hbox{\ and\ }r\ge 0\}$. $\R_1$ is a sub-semigroup of $\R$. 

\Def 3.1. Let $\mr \in \R_1$. A sequence $(X_t, t\ge 0)$ of \bzp-valued random variables is said to be an \inar1 process if for any $t \ge 0$,
$$
X_t=\mr \odot \Xtone + \eps_t,\leqno(3.1)
$$
where $(\eps_t, t\ge 1)$ is an iid sequence of \bzp-valued random variables that is assumed independent of the $Z$ variables that define the operator $\odot$ in (2.12). $\epst$ is called the innovation sequence of the \inar1 process.

The action of $\odot$ on $\Xtone$ in (3.1) is performed independently for each $t$. More precisely, we assume the existence of an array $(Z_{i,t},\ i\ge 0,\ t\ge 0)$ of iid \bzp-valued random variables, independent of $\epst$, such that the array's common pgf is $\psimr(s)$ and
$$
\mr \odot \Xtone=\sum_{i=1}^{\Xtone}Z_{i,t-1}.\leqno(3.2)
$$
These assumptions clearly make the model (3.1) a Markov chain.

In the remainder of this section $\mu_\eps$, $\sig_\eps^2$ (either or both could be infinite) and $\phi_\eps(s)$ will denote the marginal common mean, variance and pgf of the innovation sequence $\epst$ in (3.1).
  
We list several distributional properties of the \inar1 process (3.1). The proofs follow from Proposition 2.7 (see also Aly and Bouzar (1994a)).

\Prop 3.2. Let $\{X_t\}$ be an \inar1 process such that $E(X_t)<\infty$ and $Var(X_t)<\infty$ ($t\ge 0$), $\mu_\eps<\infty$ and $\sig_\eps^2<\infty$. For any $t\ge 1$,

\noi (i) $E(X_t|\Xtone)=m\Xtone+\mu_\eps$.

\noi (ii) $Var(X_t|\Xtone)=m(2r+1-m)\Xtone+\sig_\eps^2$.

\noi (iii) Let $1\le k\le t$. The covariance at lag $k$ of $\{X_t\}$ is 
$$
Cov(\Xtk,X_t)=m^kVar(\Xtk).
$$
\noi (iv) For any $t\ge 1$,
$$
E(X_t)=m^tE(X_0)+\mu_\eps\sum_{k=0}^{t-1}m^k
$$
\noi and
$$
Var(X_t)=m^{2t}Var(X_0)+(2r+1-m)\sum_{k=1}^t m^{2k-1}E(\Xtk)+\sig_\eps^2\sum_{k=1}^t m^{2(k-1)}.
$$

Next, we discuss the existence of stationary \inar1 processes. 

Since the \inar1 process (3.1) is a Markov chain, it is (strictly) stationary if and only if it admits a proper limit distribution (and it is started with that distribution). It is also a well known fact that \inar1 processes are branching processes with stationary immigration. As such, necessary and sufficient conditions for the stationarity of an \inar1 process are readily available. We list a few such conditions and refer to Foster and Williamson (1971) and Athreya and Ney (1972) for proofs and further details.

\Prop 3.3. Let $\Xt$ be an \inar1 process for some $\mr\in \R_1$.

\noi (i) If $0<\mu_\eps<\infty$, then $\Xt$ admits a proper limit distribution as $t\to\infty$.

\noi (ii) $\Xt$ admits a proper limit distribution as $t\to\infty$ if and only if
$$
\int_0^1{1-\phi_\eps(s)\over \psimr(s)-s}\,ds<\infty,\leqno(3.3)
$$
where $\phi_\eps(s)$ is the common pgf of the $\eps_t$'s (this result also holds for $m=1$).

\noi (iii) $\Xt$ admits a proper limit distribution as $t\to\infty$ if and only if $E(\ln^+ \eps_1)<\infty)$, where $\ln^+ a =\max(\ln a, 0)$, $a\ge 0$.

\Prop 3.4. Let $\Xt$ be a stationary \inar1 process for some $\mr\in \R_1$.

\noi (i) The marginal pgf $\phi_X(s)$ of $\Xt$ satisfies the equation
$$
\phi_X(s)=\phi_X(\psimr(s))\phi_\eps(s).\leqno(3.4)
$$ 

\noi (ii) Assuming $\mu_\eps<\infty$ and $\sig_\eps^2<\infty$ the correlation coefficient of $\Xt$ at lag $k$ is
$$
\rho(k)=m^k
$$

\noi (iii) The marginal mean and variance of $\Xt$ are 
$$
\mu_X={\mu_\eps\over 1-m} \hbox{\ and\ } \sig_X={m(2r+1-m)\mu_X+\sig_\eps^2\over 1-m^2}.
$$

\noi (iv) The joint pgf of $(X_{t-1}, X_t)$ is
$$
\phi_1(s_1,s_2)={\phi_X(s_1\psimr(s_2))\phi_X(s_2)\over \phi_X(\psimr(s_2))}.\leqno (3.5)
$$

\proof (i) and (iii) follow from Proposition 3.2. A standard pgf argument yields (i) and (iv). The details are omitted. \endproof

A simple induction argument (starting with (3.4)) shows that the marginal pgf of a stationary \inar1 process $\Xt$ satisfies
$$
\phi_X(s)=\phi_X(\psi_{\mrastn}(s))\prod_{k=0}^{n-1}\phi_\eps(\psi_{\mrastk}(s))\qquad (n\ge 1),
$$
with $\mrastk$ as in (2.7)-(2.8). We have by (2.8) that
$$
\phi_X(s)=\lim_{k\to \infty}\prod_{k=0}^\infty\phi_\eps(\psi_{\mrastk}(s)),\leqno(3.6)
$$
which implies the infinite order integer-valued moving average (INMA($\infty$)) representation of $\Xt$
$$
X_t=\sum_{k=0}^\infty \mrastk\odot \eps_{t-k}',\leqno(3.7)
$$
where $(\eps_j', j=0,\pm 1, \pm 2, \cdots)$ is a doubly-infinite sequence of iid random variables with common pgf $\phi_\eps(s)$.

Bourguignon and Wei{\ss} (2017) proved the existence of a stationary \inar1 process of type (3.1) with a \bergp marginal distribution, provided the parameters $m'$ and $r'$ satisfy the constraints $0<m'-r'<\min({r\over m},1)$ and $r'>{r\over 1-m}$ (stated here in terms of the re-parameterization (2.26)) . The marginal distribution of the innovation sequence $\{\eps_t\}$ is the convolution of a BerG distribution and a zero-modified geometric distribution  (see their Proposition 5).

We propose to enlarge the range of admissible values of $(m',r')$ in the Bourguignon and Wei{\ss} (2017) BerG model and we show that the marginal distribution of the innovation sequence can be written as the convolution of two BerG distributions. 

First, we need a basic result.

\Lemma 3.5. Let $\mr \in \R_1$ and $\mrp\in \R$. Then
$$
{\psimrp(s)\over\psimrp(\psimr(s))}=\psimrone(s)\psimrtwo(s),\leqno(3.8)
$$
where 
$$
\mrone=(r,r)+(r'-m')(m-1,m) \hbox{\ and\ } \mrtwo=(r'(1-m)-r,r').\leqno(3.9)
$$
Moreover, $\psimrone(s)$ is a pgf if and only if ${-r\over 1-m}<m'-r'\le \min({r\over m},1)$ and $\psimrtwo(s)$ is a pgf if and only if $r'\ge {r\over 1-m}$ (note $\psimrtwo(s)=1$ if $r'={r\over 1-m}$).

\proof Equations (3.8) and (3.9) are easily derived. In the second part of the lemma, the constraints on $m'$ and $r'$ are necessary and sufficient conditions for $\mrone$ and $\mrtwo$ to belong to $\R$. The details are omitted.\endproof

\Prop 3.6. Let $\mr\in \R_1$ and $\mrp\in \R$ such that 
$$
{-r\over 1-m}<m'-r'\le \min\bigl({r\over m},1\bigr) \quad \hbox{\ and\ } \quad r'\ge {r\over 1-m}.\leqno(3.10)
$$
Then there exists a stationary \inar1 process governed by (3.1) with a \bergp marginal distribution. The innovation sequence $\epst$ has a marginal distribution that is the convolution of a BerG$\mrone$ and a BerG$\mrtwo$, with $\mrone$ and $\mrtwo$ as in (3.9) (and noting BerG$\mrtwo$ is degenerate at $0$ if $r'={r\over 1-m}$).

\proof First, we note that by Lemma 3.5 the convolution BerG$\mrone$ $\star$ BerG$\mrtwo$ is well defined. Consider a probability space $(\Omega, {\cal{F}}, P)$ where are defined a random variable $X_0$ with a \bergp distribution, an array $\{Z_{i,t}\}$ of iid random variables with a \berg distribution, and a sequence $\epst$ of iid random variables with a BerG$\mrone$ $\star$ BerG$\mrtwo$ distribution. We assume $X_0$, $\epst$, $\{Z_{i,t}\}$ are mutually independent. Using (3.1) and (3.2), we obtain the \inar1 process $\Xt$. It follows by (3.4) and (3.8) that $X_t$ has a \bergp distribution for every $t\ge 1$. This insures stationarity of the process (by Proposition 3.3). \endproof

Barreto-Souza (2015) introduced the stationary \inar1 process 
$$
X_t=\a\otimes_{\N\B} \Xtone + \eps_t
$$  
where $\otimes_{\N\B}$ is the thinning operator of Ristic et al. (2009), the distribution of $X_t$ is a zero-modified geometric distribution (ZMG($\pi,\mu$)) with pgf 
$$
\varphi_{\pi,\mu}(s)={1+\pi\mu(1-s)\over 1+\mu(1-s)}\qquad   \mu >0 \hbox{\ and\ } -{1\over \mu}<\pi<1,\leqno(3.11)
$$
and $\a\in (\max(0,\pi\mu/(1+\pi\mu), \mu/ (1+\mu))$. The author shows that the distribution of $\eps_t$ is the convolution of two zero-modified geometric distributions, ZMG($\pi_i,\mu_i$), for some $\pi_i$ and $\mu_i$ satisfying the inequalities in (3.11), $i=1,2$.  

The following re-parameterization, 
$$  
\varphi_{\pi,\mu}(s)=\psimrp(s) \qquad m'=\mu(1-\pi) \hbox{\ and\ } r'=\mu,\leqno(3.12)
$$
shows that that the zero-modified distribution with pgf (3.11) can be seen as a BerG$\mrp$ distribution (note that $0<m' \le 1+r'$ by the inequalities in (3.11)). Let $\mr \in \R_1$. If we assume that 
$$
-\min(r/m,1) <\pi\mu < r/(1-m) \hbox{\ \ and \ \ } \mu >r/(1-m),
$$
then $m'$ and $r'$ in (3.12) satisfy (3.10). It follows by Proposition 3.6 that there exists a stationary \inar1 process of type (3.1) with the BerG$\mrp$ representation of the ZMG($\pi,\mu$) distribution as its marginal distribution.  

Next, we construct a stationary  \inar1 process with a compound negative binomial distribution.

\Lemma 3.7. Let $\mrp\in\R$ and $a>0$. If $0<m'\le r'$, then $[\psimrp(s)]^a$ is the pgf of a compound negative binomial distribution on {\bzp}. We denote such a distribution by CompNB$(m,r,a)$.

\proof It is easily seen that $\psimrp(s)=\psi_{m',m'}(\psi_{1, r'-m'}(s))$. Since $[\psi_{m',m'}(s)]^a$ is the pgf of a negative binomial distribution with parameters $(1/(m'+1),a)$, it ensues that $[\psimrp(s)]^a$ is the negative binomial compounding of iid random variables with common pgf $\psi_{1,r'-m'}(s)$. \endproof

Lemma 3.7 fails for $r'<m'\le r'+1$ as the following counterexample shows. Let $r'=0.2$, $m'=0.8$ and $a=1/2$. Then ${d^2\over ds^2}[\psi_{.8,.2}(s)]^{1/2}\big|_{s=0}=-0.24056$.

\Prop 3.8. Let $\mr\in \R_1$ and $a>0$. Assume $\mrp\in \R$ satisfies $0\le r'-m'<{r\over 1-m}$ and $r'\ge {r\over 1-m}$. Then there exists a stationary \inar1 process governed by (3.1) with a CompNB$(m',r',a)$. The innovation sequence $\eps_t$ has a marginal distribution that is the convolution of a CompNB$(m_1,r_1,a)$ and a CompNB$(m_2,r_2,a)$, with $\mrone$ and $\mrtwo$ as in (3.9).

\proof We have by (3.8)
$$
{[\psimrp(s)]^a\over[\psimrp(\psimr(s))]^a}=[\psimrone(s)]^a[\psimrtwo(s)]^a,
$$
where $\mrone$ and $\mrtwo$ are as in (3.9). The constraints $0\le r'-m'<{r\over 1-m}$ and $r'\ge {r\over 1-m}$ imply (3.10). Therefore, by Lemma 3.5, $\psimrone(s)$ and $\psimrtwo(s)$ are pgf's. Moreover, since $m_1-r_1=m'-r'\le 0$ and $m_2-r_2=-r'm-r\le 0$, we have by Lemma 3.7 that $[\psimrone(s)]^a$ and $[\psimrtwo(s)]^a$ are pgf's of CompNB$(m_1,r_1,a)$ and CompNB$(m_2,r_2,a)$ distributions, respectively. The argument that establishes Proposition 3.6 applies from this point on. The details are omitted.\endproof

If we restrict $\pi$ to the interval $[0,1)$, then Lemma 3.7 applies to the re-parameterized version (3.12) of the ZMG($\pi,\mu)$ distribution of (3.11), since $0<m'\le r'$. Therefore, letting $\mr\in \R_1$ and assuming $\pi\mu < r/(1-m)$ and $\mu >r/(1-m)$,  we can conclude by Proposition 3.8 that there exists a stationary \inar1 process with a CompNB$(m',r',a)$, where $m'$ and $r'$ are as in (3.12).  

Note that if $m'=r'$ and $r'\ge {r\over 1-m}$, then the CompNB$(r',r',a)$ \inar1 process in Proposition 3.8 is the stationary \inar1 process with a negative binomial (${1\over r'+1},a$) marginal distribution introduced by Aly and Bouzar (1994a). Moreover, the special case $m'=r'={r\over 1-m}$ gives rise to a time-reversible stationary \inar1 process with a negative binomial (${1-m\over 1-m+ r},a$) marginal distribution. Indeed, the joint pgf $\phi_1(s_1,s_2)$ of $X_{t-1}$ and $X_t$, shown to be by (3.5) 
$$
\phi_1(s_1,s_2)=\Bigl[1+r-r(s_1+s_2)+{r\over 1-m}(r(1-s)1(1-s_2)-ms_1s_2)\Bigr]^{-a}.
$$
is symmetric in $s_1$ and $s_2$, implying time reversibility. This property in fact characterizes this process as shown in Aly and Bouzar (1994a). We state the result and refer to their article for a proof (Proposition 5.1, therein).

\Prop 3.9. Let $\mr\in \R_1$. Let $\Xt$ be a stationary \inar1 process governed by (3.1) for some $\mr\in \R_1$. Assume $X_t$ has and finite mean and variance. Then $\Xt$ is time reversible if and only if its marginal distribution is negative binomial (${1-m\over 1-m+ r},a$) for some $a>0$.

\vskip .5 cm

\goodbreak\noi {\bf 4. Monotonicity}

Let $M$ and $R$ be independent random variables such that $M$ has the power distribution on $(0,1)$ with probability density function (pdf) $f_M(x)=\a x^{\a-1}$, $\a>0$, and $R$ has an exponential distribution with mean $\th>0$ and pdf $f_R(x)={1\over \th}e^{-{x\over \th}}$, $x>0$.

\Def 4.1. A \bzp-valued random variable $X$ (or its distribution) is said to be $\MRSq$-monotone if
$$
X \dover \MR\odot W,\leqno(4.1)
$$
where $W$ is a \bzp-valued random variable independent of $\MR$.

We recall (Steutel, 1988) that a \bzp-valued random variable $X$ is $\a$-monotone, $\a>0$  if
$$
X\dover (M,0)\odot W\dover M \otimes_\B W,\leqno(4.2)
$$
where $M$ is as in Definition 4.1, $W$ is a \bzp-valued random variable independent of $M$, and $\otimes_\B$ is the binomial thinning operator.

We recall two useful characterizations of $\a$-monotonicity.

\Prop 4.2. Let $(q_n, n\ge 0)$ be a pmf and $\a>0$. The following assertions are equivalent.

\noi (i) $\{q_n\}$ is $\a$-monotone.    

\noi (ii) The pgf $Q(z)$ of $\{q_n\}$ admits the representations
$$
Q(s)=\a\int_0^1 G(1-m+ms) m^{\a-1}\,dm=\a(1-s)^{-\a}\int_s^1 (1-w)^{\a-1}G(w)\,dw\leqno(4.3)
$$
for some pgf $G(s)$.

\noi (iii) For every $n\ge 0$,
$$
(n+\a)q_n\ge (n+1)q_{n+1}.\leqno (4.4)
$$

We extend Proposition 4.2 to $\MRSq$-monotonicity.

\Prop 4.3. Let $X$ be a \bzp-valued random variable with pmf $\{p_n\}$. The following assertions are equivalent.

\noi (i)  $X$ is $\MRSq$-monotone, where $M$ and $R$ are as in Definiton 4.1.

\noi (ii) The pgf $\phi(s)$ of $X$ admits the representation
$$
\phi(s)=\th^{-1}e^{1\over \th(1-s)}\int_s^1 (1-w)^{-2}e^{-{1\over \th (1-w)}}Q(w)\,dw,\leqno(4.5)
$$
where $Q(s)$ is the pgf of an $\a$-monotone distribution on {\bzp} (cf. (4.2) and (4.3)).

\noi (iii) Let $q_n=(2\th n+1)p_n - \th\bigl((n+1)p_{n+1}+(n-1)p_{n-1}\bigr)$, $n\ge 0$ (and $p_{-1}=0$). Then for every $n\ge 0$
$$
q_n\ge 0 \quad \hbox{and} \quad  (n+\a)q_n\ge(n+1)q_{n+1}.\leqno(4.6)
$$

\proof $X$ is $\MRSq$-monotone if and only if its pgf $\phi(s)$ takes the form
$$
\phi(s)=\int_0^1\Bigl(\int_0^\infty G(\psimr(s)){1\over \th}e^{-r/\th}\,dr\Bigl)\a m^{\a-1}\,dm,\leqno(4.7)
$$
where $G(s)$ is the pgf of $W$ in (4.1). The change of variable $r={1\over 1-w}-{1\over 1-s}$ in the inner integral in (4.7), along with a change of the order of integration, yield the equivalent representation 
$$
\phi(s)=\th^{-1}e^{1\over \th(1-s)}\int_s^1 \Bigl(\int_0^1 G(1-m+mw)\,\a m^{\a-1}\,dm\Bigr)(1-w)^{-2}e^{-{1\over \th (1-w)}}\,dw,\leqno(4.8)
$$
By Proposition 4.2 (first equation in (4.3)), $Q(w)=\int_0^1 G(1-m+mw)\,\a m^{\a-1}\,dm$ is the pgf of an $\a$-monotone distribution. We have thus shown (i)$\Leftrightarrow$(ii). Assume (ii) holds. Differentiating (4.5) leads to 
$$
Q(s)=\phi(s)-\th(1-s)^2\phi'(s).\leqno(4.9)
$$ 
Denoting by $\{q_n\}$ the pmf of $Q(s)$, we deduce from the power series version of (4.9) that $q_n=p_n - \th\bigl((n+1)p_{n+1}-2np_n+(n-1)p_{n-1}\bigr)$. The first part of (4.6) holds trivially and the second part follows from the fact that $Q(s)$ is $\a$-monotone and from Proposition 4.2. Thus (ii) $\Rightarrow$ (iii). We now assume that (iii) holds. Denote $d_n=np_n-(n-1)p_{n-1}$, $n\ge 1$, and $d_0=0$. Then $q_n=p_n-\th(d_{n+1}-d_n)$. This implies that $\sum_{k=0}^kq_k=\sum_{k=0}^n p_k-\th d_{n+1}$. Since $q_k\ge 0$, $\lim_{n\to \infty}\sum_{k=0}^n q_k\le \infty$. This in turn implies that $\lim_{n\to \infty} |d_{n+1}|\le \infty$. Noting that $d_{n+1}=n(p_{n+1}-p_n)+p_{n+1}$, neither $\lim_{n\to \infty} |d_{n+1}| >0$ nor $\lim_{n\to \infty} |d_{n+1}|= \infty$ can hold as that would contradict the fact that $\sum_{n=0}^\infty  |p_{n+1}-p_n|<\infty$. Therefore, $\lim_{n\to \infty} |d_{n+1}|=0$. We conclude that $\{q_n\}$ is a pmf and that it is $\a$-monotone, by the second part of (4.6). The pgf $Q(s)$ of $\{q_n\}$ must satisfy (4.9) (by definition). Solving (4.9) for $\phi(s)$ leads to (4.5). Thus (iii) $\Rightarrow$ (ii). \endproof 

One can define a notion of marginal monotonicity. 

A \bzp-valued random variable is said to be $[M,r]$-monotone if $X=(M,r)\odot W$, where $r\ge 0$ and $M$ and $W$ are as in Definition 4.1 (for some $\a>0$). The pgf $\phi(s)$ of $X$ takes the form $\phi(s)=\int_0^1 G(\psimr(s))\,\a m^{a-1}\,dm$ for some pgf $G(s)$. The change of variable (for $m$) $w=\psimr(s)$ shows that $X$ is $[M,r]$-monotone if and only if
$$
\phi(s)=\a(1-\psione(s))^{-\a}\int_{\psione(s)}^1(1-w)^{\a-1}G(w)\, dw.\leqno(4.10)
$$
We note that $\phi(s)=\phi_1(\psione(s))$, where $\phi_1(s)$ is the pgf of an $\a$-monotone distribution. $[M,r]$-monotonicity is equivalent to the $\otimes_\rho$-monotonicity of Jazi and Alamatsaz (2012).

Switching the roles of $M$  and $R$, we say that a \bzp-valued random variable is $[m,R]$-monotone if $X=(m,R)\odot W$, where $0<m\le 1$ and $R$ and $W$ are as in Definition 4.1 (for some $\th >0$). The pgf $\phi(s)$ of $X$ takes the form $\phi(s)=\int_0^\infty G(\psimr(s)){1\over \th}e^{-r/\th}\,dr$ where $G(s)$ is the pgf of $W$. Using the same change of variable as in the proof of Proposition 4.3 ((i)$\Leftrightarrow$(ii)), one can show that $X$ is $[m,R]$-monotone if and only if
$$
\phi(s)=\th^{-1}e^{{1\over \th(1-s)}}\int_s^1 G(1-m+mw)(1-w)^{-2}e^{-{1\over \th (1-w)}}\, dw.\leqno(4.11)
$$

We note that $[1,R]$-monotonicity is equivalent to the $[\Psi_1;{1\over\th}]$-monotonicity introduced by Aly and Bouzar (2002). The latter is based on the continuous semigroup of pgf's $\Psi_1=(\psione(s), r\ge 0)$ (see (2.16) and the discussion preceding it).

\Cor 4.4. Let $M$ and $R$ be as in Definition 4.1 for some $\a\ \th>0$. If $X$ is an $\a$-monotone \bzp-valued random variable, then for every $m\in (0,1]$, $(m,R)\odot X$ is $\MRSq$-monotone.

\proof Let $G(s)$ be the pgf of $X$. The pgf $\phi(s)$ of $(m,R)\odot X$ satisfies (4.11). By Proposition 4.2 applied to $G(s)$, there exists a pgf $Q(s)$ such that
$$
G(1-m+ms)=\int_0^1 Q(1-pm+pms)\,\a p^{\a-1}\,dp=\int_0^1 Q_m(1-p+ps)\,\a p^{\a-1}\,dp,
$$
where $Q_m(s)=Q(1-m+ms)$ is a pgf. Therefore, $G(1-m+ms)$ is $\a$-monotone (by appealing again to Proposition 4.2). We conclude that $\phi(s)$ admits the representation (4.5).\endproof

We note that the proof of Corollary 4.4 implies 
$$
(m,R)\odot (M\otimes_\B W)\dover \MR\odot (m\otimes_\B W)
$$
for any \bzp-valued random variable $W$.

We now address the question of monotonicity of the convolution of $\MRSq$-monotone distributions.

\Prop 4.5. Let $(M_i,R_i)$ be as in Definition 4.1 for some $\a_i, \th_i >0$, $i=1,2$. The convolution of an $[M_1,R_1]$-monotone distribution and an $[M_2,R_2]$-monotone distribution is $\MRSq$-monotone, where $M$ has the power distribution on $(0,1)$ with parameter $\a=\a_1+\a_2+1/\th_1+1/\th_2$ and $R$ has an exponential distribution with mean $\th={\th_1\th_2\over \th_1+\th_2}$.

\proof Let $\phi_i(s)$ be the pgf of the $[M_i,R_i]$-monotone distribution, $i=1,2$. Then by (4.5)
$$
H_i(s)=\th e^{-{1\over \th(1-s)}}\phi_i(s)=\int_s^1 (1-w)^{-2}e^{-{1\over \th (1-w)}}Q_i(w)\,dw,\leqno(4.12)
$$
where $Q_i(s)$ is the pgf of an $\a_i$-monotone distribution, $i=1,2$. Straightforward calculations show that
$$
{d\over ds}[H_1(s)H_2(s)]=-(1-s)^{-2}e^{{1\over 1-s}\bigl({1\over \th1}+{1\over \th_2}\bigr)}\bigl[\th_2Q_1(s)\phi_2(s)+\th_1Q_2(s)\phi_1(s)\bigr],
$$
which implies that (note $H_i(1)=0$)
$$
H_1(s)H_2(s)=\int_s^1 1-w)^{-2}e^{{1\over 1-w}\bigl({1\over \th1}+{1\over \th_2}\bigr)}\bigl[\th_2Q_1(w)\phi_2(w)+\th_1Q_2(w)\phi_1(w)\bigr]\,dw,
$$
which in turn implies (see (4.12))
$$
\phi_1(s)\phi_2(s)=(1/ \th_1+1/ \th_2)e^{{1\over 1-s}\bigl({1\over \th1}+{1\over \th_2}\bigr)}\int_s^1 1-w)^{-2}e^{{1\over 1-w}\bigl({1\over \th1}+{1\over \th_2}\bigr)}Q(w)\,dw,\leqno(4.13)
$$
where 
$$
Q(s)={\th_2\over \th_1+\th_2} Q_1(s)\phi_2(s)+{\th_1\over \th_1+\th_2}Q_2(s)\phi_1(s) \leqno(4.14)
$$
is the pgf of a two-point mixture of two distributions on {\bzp} with respective pgf's $Q_1(s)\phi_2(s)$ and $Q_2(s)\phi_1(s)$. Claim: $Q_1(s)\phi_2(s)$ is the pgf of an $[\a_1+1/\th_2]$-monotone distribution. We denote by $\{q_n\}$ (resp. $\{p_n\}$) the pmf with pgf $Q_1(s)$ (resp. $\phi_2(s)$). Let $\{(p\star q)_n\}$ be the convolution of $\{q_n\}$ and $\{p_n\}$. We have
$$
\eqalign{
(n+1)(p\star q)_{n+1}&=(n+1)\sum_{i=0}^{n+1}p_i q_{n+1-i}\cr
&= \sum_{i=0}^n (n+1-i)p_i q_{n+1-i} +\sum_{i=0}^n (i+1)p_{i+1} q_{n-i}.\cr
}
$$
Since $\{q_n\}$ is $\a_1$-monotone, we have $(n+1-i)q_{n+1-i}\le (n-i+\a_1)q_{n-i}$ (by Proposition 4.2). Therefore,
$$
(n+1)(p\star q)_{n+1}\le (n+\a_1)(p\star q)_n +\sum_{i=0}^n ((i+1)p_{i+1}-ip_i) q_{n-i}.\leqno(4.15)
$$
Since $\{p_n\}$ is $[M_2,R_2]$-monotone, we have by Proposition 4.3 ((i)$\Leftrightarrow$(iii)) 
$$
\th_2((i+1)p_{i+1}+(i-1)p_{i-1})\le (1+2\th_2 i)p_i,
$$ 
which implies that $(i+1)p_{i+1}\le \bigl({1\over \th_2}+2i\bigr)p_i$. Therefore,
$$
\sum_{i=0}^n ((i+1)p_{i+1}-ip_i) q_{n-i}\le \sum_{i=0}^n({1\over \th2}+i)p_i q_{n-i}\le (n+{1\over \th_2})(p\star q)_n.
$$
It follows from (4.15) that
$$(
n+1)(p\star q)_{n+1}\le \bigl(n+\a_1+{1\over \th_2}\bigr)(p\star q)_n,
$$
from which we conclude $\{(p\star q)_n\}$ is $[\a_1+1/\th_2]$-monotone, thus proving the claim. Using the exact same argument, one can show that $Q_2(s)\phi_1(s)$ is the pgf of an $[\a_2+{1\over \th_1}]$-monotone distribution. Since $a$-monotonicity implies $b$-monotonicity if $0<a<b$, it ensues that $Q(s)$ of (4.14) is the pgf of a two-point mixture of $[\a_1+\a_2+1/\th_1+1/\th_2]$-monotone distributions, which trivially implies the said two-point mixture is itself $[\a_1+\a_2+1/\th_1+1/\th_2]$-monotone. We conclude by (4.13) and Proposition 4.3 that $X_1+X_2$ is $[M',R']$-monotone, where $M'$ and $R'$ are independent random variables, $M'$ has the power distribution with parameter $\a=\a_1+\a_2+1/\th_1+1/\th_2$ and $R'$ has an exponnetial distribution with mean $\th=\th_1\th_2/(\th_1+\th_2)$.\endproof

Using the pgf argument in the first part of the proof of Proposition 4.5, along with (4.10) and (4.11), one can show the following holds true (cf. Proposition 4.5 for the notation).

(i) Let $r >0$. The convolution of an $[M_1,r]$-monotone distribution and an $[M_2,r]$-monotone distribution is $[M,r]$-monotone, where $M$ has a power distribution on $(0,1)$ with parameter $\a_1+\a_2$. This result is due to Jazi and Alamatsaz (2012).

(ii) Let $0<m<1$. The convolution of an $[m,R_1]$-monotone distribution and an $[m,R_2]$-monotone distribution is $[m,R]$-monotone, where $R$ has an exponential distribution with parameter $\th={\th_1\th_2\over \th_1+\th_2}$.

\bibs

\ref Alamatsaz, M.H. (1993). On discrete $\a$-unimodal distributions. Statist. Neerlandica, {\bf 47}, 245--252.

\ref Al-Osh, M. A. and Aly, E.-E.A.A. (1992) First order autoregressive time series with negative binomial and geometric marginals. Comm. Statist. Theory Methods {\bf 21}, 2483--2492.

\ref Aly, E.E. and Bouzar, N. (1994a). Explicit stationary distributions for some Galton-Watson processes with immigration. Commun. Statist.--Stochastic Models {\bf 10}, 499--517.

\ref Aly, E.-E.A.A. and Bouzar, N. (1994b). On some integer-valued autoregressive moving average models. J. of Multivariate Analysis {\bf 50}, 132--151.

\ref Aly, E.-E.A.A. and Bouzar, N. (2002). A notion of $\a$-monotonicity with generalized multiplications. Ann. Inst. Stat. Math. {\bf 54}, 125--137.

\ref Athreya, K.B. and Ney, P.E. (1972). {\it Branching Processes}. Springer-Verlag Berlin Heidelberg.

\ref Barreto-Souza, W. (2015). Zero-modified geometric \inar1 process for modelling count time series with deflation or inflation of zeros. J. Time Ser. Anal. {\bf 36}, 839--852.

\ref Bourguignon M. and Wei{\ss}, C.H. (2017). An INAR(1) process for modeling count time series with equidispersion, underdispersion and overdispersion. TEST, https:\ //doi.org/10.1007/s11749-017-0536-4. 

\ref Borges, P., Molinares, F.F., and Bourguignon, M. (2016). A geometric time series model with inflated-parameter Bernoulli counting series. Statist. Probab. Lett. {\bf 119}, 264--272.

\ref Borges, P., Bourguignon, M. and Molinares, F.F. (2017). A generalized NGINAR(1) process with inflated-parameter geometric counting series. Aust. N.Z. J. Stat. {\bf 59}, 137--150.

\ref Foster, J.H. and Williamson, J.A. (1971). Limit theorems for the Galton-Watson process with time-dependent immigration. Z. Wahrsch. Verw. Gebiete {\bf 20}, 227--235.

\ref Harris, T.E. (1963). {\it The Theory of Branching Processes}. Springer-Verlag Berlin Heidelberg.
 
\ref van Harn, K., Steutel, F.W., and Vervaat, W. (1982). Self-decomposable discrete distributions and branching processes. Z. Wahrscheinlichkeitstheor. Verw. Gebiete {\bf 61,} 97--118.

\ref Jazi, M.A. and Alamatsaz, M.H. (2012). Two new thinning operators and their applications. Glob. J. Pure Appl. Math. {\bf 8}, 13--28.

\ref McKenzie, E. (2003). Discrete variate time series. Shanbhag, D. N. (ed.) et al., Stochastic processes: Modelling and simulation. Amsterdam: North-Holland. Handb. Stat. {\bf 21}, 573-606. 

\ref Risti\'c, M.M., Bakouch, H.S., and Nasti\'c, A.S. (2009). A new geometric first-order integer-valued autoregressive (NGINAR(1)) process. J. Stat. Plann. Inf. {\bf 139}, 2218-2226.

\ref Scotto, M.G., Wei{\ss}, C.H., and Gouveia, S. (2015). Thinning-based models in the analysis of integer-valued time series: a review. Statistical Modelling, {\bf 15}15, 590--618.

\ref Steutel, F.W. (1988). Note on discrete $\a$-unimodality, Statist. Neerlandica {\bf 48}, 137--140.

\ref Steutel, F.W. and van Harn, K. (1979). Discrete analogues of self-decomposability and stability, Annals of Probability {\bf 7}, 893--899.

\ref Wei{\ss}, C. H. (2008). Thinning operations for modeling time series of counts—a survey. AStA Adv. Stat. Anal. {\bf 92} (2008), 319--341.

\ref Zhu, R. and Joe, H. (2003). A new type of discrete self-decomposability and its application to continuous-time Markov processes for modeling count data time series. Stochastic Models {\bf 19}, 235--254.

\end